\documentclass{amsart}
\usepackage{graphicx}
\usepackage{amsfonts}
\usepackage{amssymb}
\usepackage{mathrsfs}
\usepackage{hyperref}
\usepackage{xypic}
\newtheorem{thm}{Theorem}[section]
\newtheorem{cor}[thm]{Corollary}
\newtheorem{lem}[thm]{Lemma}

\theoremstyle{definition}

\theoremstyle{remark}
\newtheorem{rem}[thm]{Remark}
\numberwithin{equation}{section}

\begin{document}
\title[Yushan Jiang]{Identifiability and Global Stability Analysis on Some Partial Differential Algebraic System}
\author[YS Jiang]{Yushan Jiang}
\author[QL Zhang]{Qingling Zhang}
\address[YS Jiang,QL Zhang]{Institute of System Science Northeastern University, Shenyang,China.}%
\email{jys@neuq.edu.cn,qlzhang@mail.neu.edu.cn}%
\begin{abstract}
  We analysis some singular partial differential equations systems(PDAEs) with boundary conditions in high dimension bounded domain with sufficiently smooth boundary. With the eigenvalue theory of PDE the systems initially is formulated as an infinite-dimensional singular systems. The state space description of the system is built according to the spectrum structure and convergence analysis of the PDAEs. Some global stability results are provided. The applicability of the proposed approach is evaluated in numerical simulations on some wetland conservation system with social behaviour.
\end{abstract}
\maketitle
\section{Introduction}
  Singular systems have abilities in representing dynamical systems in the areas of electrical circuits and multibody systems, chemical engineering and economic systems, mechanical structure and biological systems. Singular systems are also referred to as  descriptor systems, differential-algebraic systems(DAEs), implicit systems or generalized state-space systems\cite{Dai1989}. A large number of fundamental ideas and results based on state-space systems have been successfully extended to singular systems, such as controllability observability, pole assignment, stability and stabilization \cite{Campbell1982, Lewis1986, Dai1989, Zhang1997, Riaza2007, Zhang2012} etc.
  In recent years, partial differential algebraic equations (PDAEs) become an independent field of research which is gaining in importance and becoming of increasing interest for both applications and mathematical theory.The PDAEs research areas can be classified into three groups:
 \begin{itemize}
   \item  Index analysis and solvability of PDAEs.(see \cite{Ali20120002,Ali20104666,CAMPBELL1999,Rang2005437,Chudej2005,Lennart2013,Lucht2005129002402})
   \item Innovative and improved numerical methods to solve PDAEs.(see \cite{Debrabant2005213,Vuong2014115,Bartel201414})
   \item Control problems and optimization described by PDAEs. (see\cite{HuaiNing20111172,Biao20126198365,Moghadam20130024, Daafouz201492,Tang20110002,Clever20120003,Jadachowski201585})
 \end{itemize}
  During the last few years there has been a tremendous amount of activity on PDAEs. Most practical industrial processes inherently distributed in space and time involves the use of PDAEs. Examples include integrated circuits \cite{Ali20104666}, Chemical Reactor \cite{VuTienDung2007201} and population systems \cite{Yushan2008}. Compared with DAEs there exist few results on PDAEs reported in literature.

  The structure analysis and solvability on PDAEs can be traced back to \cite{CAMPBEL1995TJ35800002}. With the same method of line (MOL) there existed different spatial indexes between different numerical methods. In order to relate properties of the PDAEs to those of the resulting DAEs it is necessary to have a concept of the index of a possibly constrained PDAEs. Firstly in \cite{CAMPBELL1999} the perturbation index was defined on a class of liner time invariant PDAEs. On infinite domain the wave solutions of PDAEs was also considered\cite{Campbell1997}. In\cite{WIESLAW1997, Marszalek20027667152} there was a systematically analysis about three different type of indexes(model, perturbation and algebraic). In \cite{Lucht2002317036213} a consistent representation of the solution of an initial boundary value problem for PDAEs was proposed. The index involved in the problem is characterized by means of the Fourier and Laplace transformations. The index jump was also discussed. The differentiation index \cite{MARTINSON20002295} of some general nonlinear PDAEs on hyper-plane domains is a generalization of the differentiation index of DAEs. The differentiation index provide a way to determine a Cauchy data on domain surface which must be consistent with the PDAEs. For the first time a perturbation index for a singular PDE of mixed parabolic-hyperbolic type was computed by \cite{Rang2005437}. \cite{Riaza2007} gave some index of PDAEs with the sequence of matrices method. It is an effective method for the finite systems. But the question is the matrices sequence built in PDAEs are unbounded operators. However, it generalizes the Kronecker index in a rather functional analytic manner. Other researchers also investigated the index structure of PDAEs systems in the area of coupling nonlinear PDAEs systems \cite{Ali20120002}, mixed index \cite{Lennart2013}, index determination algorithm \cite{Lamour2013}.

  Despite the complexities analysis of index on PDAEs, research in the area of control problems for PDAEs is relatively scarce. This approach was first applied in \cite{Zhang2004380} to design energy based controller of a coupled wave-heat equation systems. In \cite{Tang20110002}stabilizing a coupled PDE-ODE systems with interaction at the interface with boundary control was considered. \cite{Daafouz201492}designed a nonlinear saturating control law using a Lyapunov function for the averaged model of the switched power converter system. Based on the infinite-dimensional state-space representation theory \cite{Moghadam20130024} addressed the linear quadratic regulator control of the PDAEs. The optimal control problem is treated using operator Riccati equation approach. Thought the previous methods\cite{Clever20120003,Lamour2013} are derived from DAEs theory. Other optimal control problem\cite{Biao20126198365, Reis2014008, Jadachowski201585} are considered with parabolic type PDAEs.

  However in these coupled PDE-ODE systems, parabolic distributed PDEs systems or hyperbolic systems above, the time state variables matrix is reversible and the spatial variables performs in one-dimensional interval. Some singular systems like parabolic-elliptic partial differential equation can not be direct application of the above theoretical study. Our study derived from the article \cite{Lucht2005129002402} in which the search for series solutions of PDAEs with method on DAEs. Motivated by the technique in \cite{HuaiNing20111172}, we consider some singular partial differential equations systems with boundary conditions in high dimension bounded domain. The present work focuses on the development of an generalization stability method for a class of PDAEs with singular time derivative coefficient matrix. In such a systems the spatial variables act on the bounded high dimension domain.

  The organization of the study  is as follows: The problem statement for some parabolic PDAEs is given in section 2. In section 3 the original PDAEs is described as an infinite-dimensional singular systems. Then the state variables expression is built with Kronecker-Weierstrass form. And the spectrum analysis is given to show the analytical solution of PDAEs is convergence. The dynamical stability properties are analysed with Lyapunov method. Finally, in section 4 as an  application we build some PDAEs model on some wetland conservation system with social behaviours. And the stability property of the ecosystem is given by our developed PDAEs theory.

  Notations: $\mathbb{R},\mathbb{R}^{n}$ and $\mathbb{R}^{n\times m}$ are the set of real numbers,  the $n-$dimensional Euclidean space and the set of all real $n\times m$ matrices respectively. $\|\cdot\|$ denotes the Euclidean norm for vectors. For a vector $x(t)\in \mathbb{R}^{n}$,$\| x(t)\|_{\infty} =: \sup\| x(t)\|$. $x(t)\in \mathscr{L}_\infty $ if $\| x(t)\|_\infty<\infty$. For a symmetric matrix $M$,$M>(<)0$ means that it is positive (negative) definite. $I$ is the identity matrix.The superscript $T$ is used for the transpose. Matrices, if not explicitly stated, are assumed to have compatible dimensions. For the convenience, we define the following Hilbert space:
  $$\mathscr{H}\triangleq \{x:\Omega\times[0,+\infty)\rightarrow \mathbb{R}^n \ \textrm{and} \ \|x\|_{2}<\infty\}$$
  with inner product and $ \mathscr{H}_{2}$-norm  respectively defined by
\begin{eqnarray}
\nonumber
\langle x_1,x_2\rangle \triangleq&
\{\int_{\Omega}x^{T}_1x_2\textrm{d}\textbf{z}\},
\|x\|_{2} \triangleq
\{\int_{\Omega}\|x\|^2\textrm{d}\textbf{z}\}^{1/2}.
\end{eqnarray}
\section{Description of PDAEs}
  We consider the linear partial differential algebra equations system (PDAEs) in $d-$dimensional bounded spatial domain with a state space description of the form
    \begin{eqnarray}
    \label{PDAE01}
    &&E\frac{\partial x}{\partial t}=D\Delta{x}+A{x}+B{u}, \\
    \label{PDAE02}
    && {y}(t)=\int_{\Omega}C{x}\textrm{d}\textbf{z}+{v}.
    \end{eqnarray}
  subject to the boundary conditions(BCs)
    \begin{equation}\label{PDAE03}
    p_\iota x_\iota(t,\textbf{z})+q_\iota\frac{\partial x_\iota(t,\textbf{z})}{\partial \textbf{n}}=0, \textbf{z}\in\partial\Omega,\iota\in\vartheta_{BC}
    \end{equation}
  and the initial conditions(ICs)
    \begin{equation}\label{PDAE04}
    x_\kappa(0,\textbf{z})=x^0_\kappa(\textbf{z}),\textbf{z}\in\Omega,\kappa\in\vartheta_{IC}
    \end{equation}
  where
    $E,A,D\in \mathbb{R}^{n\times n}, D\geq 0$,
    $B \in \mathbb{R}^{n\times n_u}$,
    ${x}=[x_1(t,\textbf{z}) \cdots x_n(t,\textbf{z})]^{T}\in \mathbb{R}^{n}$
  is the vector of state variables,
    $\Omega \subset\mathbb{R}^d $
  is the spatial domain of the process,
    $ \textbf{z}=(z_1,\cdots,z_d)\in\Omega$
  is the spatial coordinate,
    $t\geq 0$
  is the time,
    $\Delta=\Sigma_{i=1}^{d}\partial^{2}_{z_i}$
  is the Laplace operator,
    ${u} (t,\textbf{z})\in \mathbb{R}^{n_u}$
  is the manipulated input vector function,
    ${y}(t)\in\mathbb{R}^{n_y}$
  is the measure output,
    ${v}(t)$
  is the measurement disturbance.
    $p_\iota,q_\iota$
  are given constants,
    $\textbf{n}$
  is the outward normal vector of ${x}$ on $\partial\Omega$,
    $\vartheta_{BC}\subseteq\{1,\cdots,n\}$,
    $\vartheta_{IC}\subseteq \{1,\cdots,n\}$,
    $x^0_\kappa(\textbf{z})$
  is the spatial initial state function.
\begin{rem}
  For these systems we show that at least one of the matrices $E,D$ is singular. The special case $D=0$ leads to DAEs. Therefore in this study we assume that $D$ is not a zero matrix. Another special case $|E|\neq 0$ leads to the PDE-ODE systems which some researchers investigated\cite{Daafouz201492}. In contrast to the problems with nonsingular matrices $E$ and $D$, the BCs (\ref{PDAE03}) and ICs (\ref{PDAE04}) have to fullfill certain consistency conditions. Thus, for given PDAEs (\ref{PDAE01}) the subsets
    $\vartheta_{BC},\vartheta_{IC}$
  are determined by consistency conditions(see\cite{Lucht2002317036213}). In general, the BCs  for $\iota\notin \vartheta_{BC}$ and the ICs for $\kappa\notin \vartheta_{IC}$ must be determined with the help of the PDAEs(\ref{PDAE01}). With different $p_\iota,q_\iota$ values the BCs can be three boundary types: Dirichlet, Neumann and Robin boundary conditions. We assume that the BCs are homogeneous for simplicity, thought our propose is applicable to the general inhomogeneous boundary conditions. Moreover, we assume that the process in (\ref{PDAE01})-(\ref{PDAE04}) evolves on a compact set, i.e.,
    $({x}(t,\textbf{z}),{u}(t))$ $\in X\times U$ for all $t\geq 0$,
  where
    $X\times U \subset \mathbb{R}^{n} \times \mathbb{R}^{n_u}$
  is a compact set containing the origin.
\end{rem}
\section{Decomposition of PDAEs and spectrum analysis}
\subsection{Decomposition of PDAEs and infinite singular systems}
  In this section,the PDAEs (\ref{PDAE01})-(\ref{PDAE04}) will be rebuilt as a large infinite dimensional systems. Consider the linear elliptic operator
    \begin{eqnarray}\label{Stru_liu01}
    \mathscr{D}{x}=D\Delta {x}+A{x}
    \end{eqnarray}
  in $\Omega$ with homogeneous BCs (\ref{PDAE03}). We denote $\sigma( \mathscr{D})$ be the spectrum of $ \mathscr{D}$ that is  the set of $\lambda\in \mathbb{C}$ for which $( \mathscr{D}-\lambda I)$ is not invertible. More specifically, for the operator $\mathscr{D}$, the eigenvalue problem is defined as
    \begin{equation}\label{Sect3_01}
    \mathscr{D}\phi_j(\textbf{z})=\lambda_j\phi_j(\textbf{z}), (j=1,\cdots, \infty),\phi_j(\textbf{z})\in \mathscr{H}_{2,\Omega}\ \textrm{ with}\ \textrm{BCs}\ (\ref{PDAE03}).
    \end{equation}
  where $\lambda_j\in \mathbb{R}$ denotes the $j$th eigenvalue and $\phi_j(\textbf{z})$ denotes the corresponding orthonormal eigenfunction. Obvioursly, the eigenfunctions
    $\{\phi_j(\textbf{z})\}_{j=1}^{\infty}$
  form an orthonormal basis for eigenvalue problem (\ref{Sect3_01}). Furthermore in the next subsection we show that $\sigma( \mathscr{D})$ has dscrete spectrum consisting only of real eignevalues with at most a finite mumber of positive eigenvalues which is a generalization of \cite{HuaiNing20111172}. Applying for PDE theory, we can obtain the following infinite-dimensional differential algebric equations systems(IDAEs):
    \begin{eqnarray}
    \label{DAE01}
    &&E\dot{X_j}=\lambda_{j}X_j+BU_j, j=1,\cdots,\infty \\
    \label{DAE02}
    &&Y_j=C_jX_j+V_j
    \end{eqnarray}
  with the initial condition
    \begin{equation}\label{DAE03}
    X_j(0)=\langle{x}_0(\textbf{z}),\phi_j(\textbf{z})\rangle
    \end{equation}
  where
  \begin{equation}\label{DAE06}
    X_j=\langle{x},\phi_j \rangle,
    U_j=\langle{u},\phi_j(\textbf{z})\rangle,
    V_j=\langle{v},\phi_j(\textbf{z})\rangle,
    C_j=C\int_{\Omega} \phi_j(\textbf{z})\textrm{d}\textbf{z}.
  \end{equation}

\begin{rem}
  One can not guarantee that each system are solvable without considering the regularity of the matrix pencil set $\{(sE,\lambda_{j}I)\}_{j=1}^{\infty}$. Thus in this study we assume that for each $j$ the matrix pencil with respect to the system (\ref{DAE01}) and (\ref{DAE02}) is  regular. Furthermore, all systems have the uniform differential time index\cite{Lucht2005129002402}, that is, for all $j$, the pencil $(sE,\lambda_{j}I)$ is regular and has the Riesz index $\nu_{d,j}=\nu_{d}$ (independent of $j$).
\end{rem}
\subsection{Eigenvalue estimation and property}
  For the unforced PDAEs
    \begin{eqnarray}
    \label{PDAE05}
    E\frac{\partial{x}}{\partial t}=D\Delta{x}+A{x}.
    \end{eqnarray}
  with the corresponding IDAEs
    \begin{eqnarray}
    \label{DAE04}
    E\dot{X_j}=\lambda_{j}X_j, j=1,\cdots,\infty.
    \end{eqnarray}
 For every $j$, firstly let us denote $\sigma(E,\lambda_j I)$ be the generalized eigenvalue set of systems (\ref{DAE04}) where $\lambda_j$ is defined by (\ref{Sect3_01}). We also write $\sigma(E)=\sigma(I,E)$ is the usual eigenvalue set. And for every component $x_k(t,\textbf{z})$ of state vector ${x}$, $\mu_j^k$ denotes the eigenvalue of the Strum-Liouville problem
    \begin{equation}\label{scalar SL}
    \Delta x_k(t,\textbf{z})+\mu_j^kx_k(t,\textbf{z})=0 \textrm{ subject}\ \textrm{to}\ (\ref{PDAE03}).
    \end{equation}
  The subsequent corollary concerning the eigenvalues estimation of $\lambda_j$ for (\ref{DAE04}) corresponding to systems (\ref{PDAE01})-(\ref{PDAE04}). Before that, we give two lemmas about the spectrum structural property of the scalar elliptic operator.
\begin{lem}[\cite{Smoller1994}]\label{scalar elliptic operator}
  With the strongly scalar elliptic operator $P$ ($P$ is strongly elliptic in the sense that there exists $\gamma>0$ such that for all $\xi$,
    $\sum_{i,j} p_{i,j}(\textbf{z}) \xi_i \xi_j \geq \gamma |\xi|^2$)
  the operator defined by
  $\mathscr{P} x
  \equiv Px+ax
  \equiv \sum^{d}_{i,j=1}
  \partial_{z_j}(p_{i,j}(\textbf{z})\partial_{z_i}x)
  +a(\textbf{z})x$
  in the bounded domain $\Omega\in \mathbb{R}^d$ with homogeneous boundary conditions of the form
    \begin{eqnarray}
    \alpha(\textbf{z})\frac{\textrm{d} x}{\textrm{d} \textbf{n}}+\beta(\textbf{z})x=0\ \ on\ \partial\Omega.
    \end{eqnarray}
  has discrete spectrum consisting only of real eigenvalues. If $a(\textbf{z})$ is bounded in $\Omega$, $\mathscr{P}$ can have at most a finite number of positive eigenvalues.
\end{lem}
\begin{lem}[The Bauer-Fike Theorem]
  Assume that $A\in C^{n,n} $ is a diagonalizable matrix, $V\in C^{n,n} $ is the non-singular eigenvector matrix such that $A=V\Lambda V^{-1}$, where
    $\Lambda=diag(\mu_1,\mu_2,\cdots,\mu_n)$
  is a diagonal matrix, Let $\lambda$ be an eigenvalue of $A+T$ then there exists
    $\mu\in \sigma(A)$
  such that
    $|\mu-\lambda|\leq ||V||_{2}||T||_{2}||V^{-1}||_{2}.$
\end{lem}
\begin{cor}\label{coll_35}
  For the eigenvalue problem (\ref{Sect3_01}) the operator $\mathscr{D}$ has discrete spectrum consisting only of real eigenvalues. And there is a finte number of positive eigenvalues,i.e.,if all eigenvalues $\lambda_j$ are ordered that $\lambda_j>\lambda_{j+1}$, then there exist a finite number $g$ so that $\lambda_{g+1}<0$ and
    $\varepsilon \triangleq \frac{|\lambda_1|}{|\lambda_{g+1}|}<1$
is a small positive number.
\end{cor}
\textbf{Proof}:
  Considering the Strum-Liouville problems
    \begin{eqnarray}
    \label{DAE05}
    \Delta {x}_k-\mu_j^{k}{x}_k=0\ \textrm{with BCs}\ (\ref{PDAE03}) (j=1,\cdots, \infty, k=1,\cdots,n)
    \end{eqnarray}
  where $\mu_j^{k}$ denotes the eigenvalue of (\ref{DAE05}).
  It is true that all $\mu_j^{k}$ satisfy Lemma \ref{scalar elliptic operator}. And from the assumption $D\geq 0$ one can get that for every given $j$ if for all $k=1,\cdots,n, \mu_j^{k}<0$ then all eigenvalue of $D\Delta{x}$ is negative defined. Therefore, it is following from Lemma 3.3 that for every $j$, $\mu^{k}_{j}$ and $\lambda_{j}$ satisfy
    $|\mu^{k}_{j}-\lambda_j| \leq c$,
  where $c$ is some given constant independed with $j$. Noticing that
    $\mu^{k}_{j}\rightarrow -\infty (j\rightarrow \infty)$,
  the conclusion holds.
\begin{rem}
  It should be pointed out that the existing one dimension result \cite{HuaiNing20111172} can not be directly generalized to high dimension case since the mathematical complexity property of the Strum-Liouville(L-S) problem on high dimension spatial space.  Additionally, the zero eigenvaule $\lambda_j=0$ yields the trivial constant state response ${x}(t,\textbf{z})=\textrm{const.}$ which is not interested in our study.
\end{rem}
  From Corollary\ref{coll_35} there exists some $p\in \mathbb{Z}^{+}$ such that the operator $\mathscr{D}$ with respect to (\ref{Sect3_01}) has discrete spectrum structure
    $$\sigma(\mathscr{D})
    =\{\lambda_1 \geq \lambda_2 \geq \cdots \geq \lambda_p \geq \lambda_{p+1} \geq \cdots \}$$
  and $\lambda_{p+1}<0$. Now for the matrix pencil $(E,\lambda_jI)$ with respect to IDAEs (\ref{DAE01}) if
    $s_j\in \sigma(E,\lambda_j I)$,$e\in\sigma(E)$
  then from the definition of eigenvalue the relationship among $\lambda_j,s_j$ and $e$ is
    $s_j=\frac{\lambda_j}{e}.$
Immediately we have the following spectrum structure theorem.
\begin{thm} \label{thm_eigenvalue anal}
  Assume $\{\lambda_j\}_{j=1}^{\infty}$ is the spectrum set of the S-L problem (\ref{Sect3_01}) with respect to systems (\ref{DAE02}). $e\in \sigma(E)$ is the eigenvalue of $E$, then there exist some positive number $p>0$ such that for every $j\in \mathbb{Z}^+$,$j>p$, the sign of the real part of $s_j$ is reversal with $e$. Consequently, $E$ is a nonnegative matrix is the necessary condition for which the matrix pencil family set
    $\{(E,\lambda_jI)\}_{j=1}^{\infty}$
  is admissible.
\end{thm}
\subsection{State and output response of PDAEs}
  In what follows we shall assume that each pencil $(E,\lambda_jI)$ with respect to (\ref{DAE01}) is regular and impulse-free. Intuitively, the Kronecker-Weierstrass equivalent form of DAEs\cite{Dai1989} can be applied to our systems. However, considering the nonnegative diffusion matrix $D$ this equivalent transformation leads to an unexpected matrix. And the spectrum analysis theorem above is no longer valid. Thus the Jordan transformation is introduced to solve the problem.

  For the IDAEs (\ref{DAE01}), there exists non-singular matrix $M$ such that $(E,\lambda_jI)$ can be transformed into the following form
 \begin{equation}\label{sect3_08}
    M^{-1}[E,\lambda_jI] M
    =\left[
      \left(
        \begin{array}{cc}
            E_1 & 0 \\
            0 & J \\
        \end{array}
      \right),
      \left(
        \lambda_j I
      \right)
    \right]
 \end{equation}
  where $J\in \mathbb{R}^{n\times n}$ is a  nilpotent matrix with index $n_0$, $E_1 \in\mathbb{R}^{r\times r}$ is a nonsingular matrix. It follows immediately from the linear singular systems theory\cite{Zhang2012} that the state response of such sytem is
    $$X_{j}(t)=[X_{j,\textrm{reg}}(t),X_{j,\textrm{nil}}(t)]^{T}$$
  where $X_{j,\textrm{reg}}(t),X_{j,\textrm{nil}}(t)$ are given by
 \begin{equation}\label{solution01}
 X_{j,\textrm{reg}}(t)
 =M\textrm{e}^{\lambda_j E_1^{-1}t}X_{j,\textrm{reg}}(0)
 +M\int^{t}_{0}\textrm{e}^{\lambda_j E_1^{-1}(t-s)}B_1U_j(s) \textrm{d}s
\end{equation}
and
\begin{equation}\label{solution02}
 X_{j,\textrm{nil}}(t)=-M\sum^{n_0-1}_{i=0}\lambda_j^{-(i+1)}J^{i} B_2U_j^{(i)}(t).
\end{equation}
$B_1,B_2$ are the components of $M*B$ corresponding to the partition of the vector $X_{j}(t)$ into the regular and nilpotent parts. Thus the PDAEs (\ref{PDAE01})-(\ref{PDAE04}) have the state response
$${x}(t,\textbf{z})=\sum_{j=1}^{\infty} X_{j}(t)\phi_j(\textbf{z})$$
and the output response
    $${y}(t)
    =\sum_{j=1}^{\infty}(C_jX_{j}(t)+V_j(t))
    $$
where $C_j,V_j(t),(j\in \mathbb{Z}^+)$ are defined by (\ref{DAE06}).

The above state and output reponse discuss provides a generalized systems theoretical approach to the PDAEs. From (\ref{solution01}) we know that if there exists an eigenvalue of $E_1$ with negative real part then the solution grows exponentially which coincides with the `explosive solution' from  a mathematical perspective.

It should be noted that the existing works in\cite{Moghadam20130024,Tang20110002}  rely on the assumptions that the derivative matrix $E$ is nonsingular. And for the orthogonal property of the eigenfunction set $\{\phi_j\}^{\infty}_{j=1}$, some dynamical properties of the PDAEs including stability, stabilizability, and detectability can be given through a direct application of the LMI technique about the generalized systems theory\cite{Zhang2012,Yang2013}. The following theorem shows the admissible of PDAEs via LMIs.
\begin{thm}[Admissible Via LMIs]\label{theorm_lyap ineqn}
The PDAEs (\ref{PDAE05}) is admissible if and only if there exists a matrix $P$ such that the first subsystem of IDAEs(\ref{DAE04}) is admissible,i.e.
\begin{displaymath}
\begin{array}{l}
  E^{T}P=P^{T}E\geq 0, \\
  \lambda_1(P^{T}+P) +P^{T}A+A^{T}P<0
\end{array}
\end{displaymath}
where $\lambda_1$ is the maximum eigenvalue of the spatial differential operator (\ref{Stru_liu01}).
\end{thm}
\textbf{Proof:}
Noticing that all eigenvalues $\lambda_j$ be ordered as $\lambda_{j}\geq\lambda_{j+1}$, it can be easily deduced from above that
\begin{displaymath}
  \lambda_j(P^{T}+P) +P^{T}A+A^{T}P<0
\end{displaymath}
holds for every $j\in \mathbb{Z}^+$. Then the desired result follows immediately by the admissible theory of the pencil $(E,A)$(see \cite{Xu2006} for detail).
\section{Stability analysis}
 In view of theorems \ref{thm_eigenvalue anal} and \ref{theorm_lyap ineqn} above, the systems stability can be determined by spectrum analysis and LMIs. This validity of the above stability analysis relies on the convergence properties about the IDAEs. In the following, we propose some exponential stability property on the PDAEs by the energy estimation theory. For the following considerations, it will be simplest to assume homogeneous Neumann boundary conditions.
\begin{lem}[\cite{Smoller1994}]\label{lem_Poincare}
Let $x\in W^{2}_{2}(\Omega)$,then if $\mu_1$ is the smallest positive eigenvalue of $-\Delta$ on $\Omega$ (with the appropriate boundary conditions) the following Poincar\'{e} inequalities hold:
\begin{equation}\label{Poincare inequalities01}
  \|\nabla x\|^2\geq \mu_1 \|x-\bar{x}\|^2, \|\Delta x\|^2\geq \mu_1\|\nabla x\|^{2}\ if\ \frac{dx}{d\textbf{n}}=0\ on \ \partial\Omega;
\end{equation}
\begin{equation}\label{Poincare inequalities02}
  \|\nabla x\|^2\geq \mu_1 \|x\|^2 \ if\ x=0\ on \ \partial\Omega,
\end{equation}
where $\bar{x}=\frac{1}{|\Omega|}\int_{\Omega}xd\textbf{z}.$
\end{lem}
\begin{thm}\label{thm_stable_analy}
Assume that ${x}(t,\textbf{z})$ is a bounded solution of (\ref{PDAE01}),(\ref{PDAE02})with homogeneous Neumann boundary conditions. Assume that $E\geq 0, D>0, \mu_1$ is the smallest positive eigenvalue of $-\Delta$ on $\Omega$, $d_1$ is the smallest positive eigenvalue of $D$ and
$$\delta=\frac{2(d_1\mu_1-\|A\|)}{\|E\|}>0.$$
Then
\begin{equation}\label{theorm_stable01}
\frac{1}{2}\int_{\Omega}\sum^{d}_{i=1}\frac{\partial {x}^{T}}{\partial z_i}E\frac{\partial {x}}{\partial z_i}\textrm{d}\textbf{z}\leq c_1e^{-\delta t}
\end{equation}
for a positive constant $c_1$, and
\begin{equation}\label{theorm_stable02}
\int_{\Omega}\|{x}(t,\textbf{z})-{x}_M(t)\|\textrm{d}\textbf{z}\leq c_2e^{-\delta t}
\end{equation}
where $c_2$ is a positive constant, ${x}_M(t)=\frac{1}{|\Omega|}\int_{\Omega}{x}(t,\textbf{z})\textrm{d}\textbf{z}$
is the spatial average function, i.e. the state variable vector ${x}(t,\textbf{z})$ generated by the PDAEs is exponentially stable and  asymptotically converge to its spatial average.
\end{thm}
\textbf{Proof:}By introducing the energy integral(lyapunov function)
\begin{equation}\label{Energy function}
  E_{L}(t)=\frac{1}{2}\int_{\Omega}\sum^{d}_{i=1}\frac{\partial {x}^{T}}{\partial z_i}E\frac{\partial {x}}{\partial z_i}\textrm{d}\textbf{z}
\end{equation}
and compute
\begin{eqnarray}
  \frac{\textrm{d}}{\textrm{d}t}E_L(t) &=&\int_{\Omega}\sum^{d}_{i=1}\frac{\partial {x}^{T}}{\partial z_i}E\frac{\partial^2 {x}}{\partial z_i\partial t}\textrm{d}\textbf{z}\\
  \nonumber
         &=&  \int_{\Omega}\sum^{d}_{i=1}\left(
      \frac{\partial {x}^{T}}{\partial z_i}\frac{\partial}{\partial z_i}( D\Delta {x}+A{x})
      \right)
      \textrm{d}\textbf{z}\\
  \nonumber
      &=&  \int_{\Omega}\sum^{d}_{i=1}\left(
      \frac{\partial {x}^{T}}{\partial z_i}D\frac{\partial \Delta{x}}{\partial z_i}
      +
      \frac{\partial {x}^{T}}{\partial z_i}A\frac{\partial {x}}{\partial z_i}
      \right)
      \textrm{d}\textbf{z}
\end{eqnarray}
Notice that the Neumann BCs are $\frac{\partial{x}}{\partial\textbf{n}}=0\ \textrm{on}\ \partial\Omega$. With the application of divergence theorem we have
 \begin{eqnarray}
  \frac{\textrm{d}}{\textrm{d}t}E_L(t) &=&
  -\int_{\Omega}\Delta{x}^{T}D\Delta{x}\textrm{d}\textbf{z}
   +\int_{\Omega}\sum^{d}_{i=1}
      \frac{\partial {x}^{T}}{\partial z_i}A\frac{\partial {x}}{\partial z_i}
       \textrm{d}\textbf{z}.
\end{eqnarray}
 The above equality can be estimated as
\begin{equation}\label{eqn_Energy08}
  \frac{\textrm{d}}{\textrm{d}t} E_L(t)\leq
  -d_1 \| \Delta{x}\|^2
  +\int_{\Omega}\sum^{d}_{i=1}
  \frac{\partial {x}^{T}}{\partial z_i}A\frac{\partial {x}}{\partial z_i}
  \textrm{d}\textbf{z}
\end{equation}
 where $d_1$ is the smallest positive eigenvalue of $D$. According to Lemma \ref{lem_Poincare}, (\ref{eqn_Energy08}) implies
\begin{equation}\label{eqn_Energy09}
  \frac{\textrm{d}}{\textrm{d}t} E_L(t)\leq
  -d_1 \mu_1 \| \nabla{x}\|^2+\|A\|\|\nabla {x}\|^2.
\end{equation}
Noticing that $- d_1\mu_1+\|A\|<0$, and takeing into account
\begin{equation}\label{Energy function}
  E_{L}(t)=
  \frac{1}{2}\int_{\Omega}
  \sum^{d}_{i=1}
  \frac{\partial {x}^{T}}{\partial z_i}E\frac{\partial {x}}{\partial z_i}
  \textrm{d}\textbf{z}
  \leq
  \frac{\|E\|}{2}\int_{\Omega}
  \sum^{d}_{i=1}
  \frac{\partial {x}^{T}}{\partial z_i}\frac{\partial {x}}{\partial z_i}
  \textrm{d}\textbf{z}
  =\frac{\|E\|}{2}\| \nabla{x}\|^2
\end{equation}
we obtain that
\begin{equation}\label{eqn_Energy10}
  \frac{\textrm{d}} {\textrm{d}t} E_L(t)\leq -\delta E_L(t)
\end{equation}
where
$$\delta=\frac{2(d_1\mu_1-\|A\|)}{\|E\|}>0.$$
By lemma \ref{lem_Poincare}, (\ref{eqn_Energy10}) implies (\ref{theorm_stable02}). Thus, under the conditions of the theorem, spatial oscillations decay exponentially, and the solution asymptotically behaves like its spatial average.
\section{Application to the Coastal Wetland Conservation System with social behaviour}
 \begin{figure}
  \includegraphics[width=0.7\textwidth]{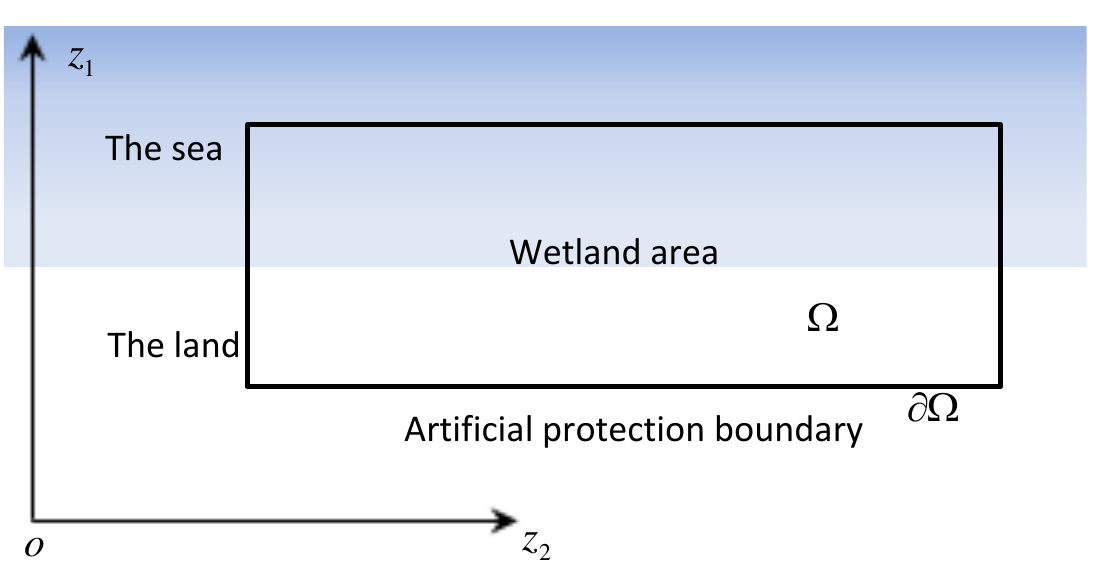}\\
  \caption{\label{ecologysys}The idealized spatial domain is a rectangular domain with the sea oriented direction $z_1$ and coast line direction $z_2$.The wetland conservation is closed with no flux boundary conditions imposed.}
\end{figure}
In this section, we illustrate, through computer simulations, the application of the theoretical development given above to some wetland conservation system with social behaviour on some plane rectangle domain (see Fig. \ref{ecologysys}).

Our model derived from \cite{Ko20139} in which a reaction-diffusion system incorporating one prey and two competing predator species under homogeneous Neumann boundary conditions was considered. In this model we choose human, birds(predators) and their food (prey) as the research objects. The spatiotemporal dynamics between predator and their prey with human activity affect in a protected environment can be described by the following PDAEs
\begin{equation}\label{prey_predator systems}
\left\{
\begin{aligned}
  &\frac{\partial x_1}{\partial t}
  =d_1\Delta x_1+r_1x_1(1-\frac{x_1}{N_1}-k_1 \frac{x_2}{N_2}-h_1x_3), \textbf{z}\in\Omega,t>0,\\
  &\frac{\partial x_2}{\partial t}
  =d_2\Delta x_2+r_2x_2(-1+k_2\frac{x_1}{N_1}-\frac{x_2}{N_2}-h_2x_3), \textbf{z}\in\Omega,t>0,\\
    &0
  =\Delta x_3+x_3, \textbf{z}\in\Omega,\\
    &\frac{\partial x_1}{\partial \textbf{n}}
  =\frac{\partial x_2}{\partial \textbf{n}}
  =\frac{\partial x_3}{\partial \textbf{n}}=0, \textbf{z}\in\partial\Omega,t>0,\\
     &x_1(0,\textbf{z})=x_1^0(\textbf{z})\geq 0,x_2(0,\textbf{z})=x_2^0(\textbf{z})\geq 0,\textbf{z}\in\Omega.\\
\end{aligned}
\right.
\end{equation}
where $x_1(t,\textbf{z}),x_2(t,\textbf{z})$ represent the population of prey (birds) and predator (fish) species at time $t>0$ and spatial position $\textbf{z}\in\Omega$ respectively; $d_1, d_2$ stand for the diffusion coefficients of prey and predator species; the prey population follows the logistic growth in the absence of predator with the intrinsic growth rate $r_1$ and the carrying capacity $N_1$; $r_2$ is the death rate of the predator with the carrying capacity $N_2$;$k_1,k_2$ represent the strength of relative effect of the interaction on the two species.

For a wetland ecosystem, the influence of human can be regard as an invasive species and not be affected by other species. Thus the influence of human activities(for example,the economic interest) on the two species are represented by $h_1x_3$ and $h_2x_3$ respectively (in the first two equations of system (\ref{prey_predator systems})) and the human population is in a state of free distribution can be described by
\begin{equation}\label{human_distr01}
\frac{\partial x_3(t,\textbf{z})}{\partial t}=\Delta x_3(t,\textbf{z})+x_3(t,\textbf{z}),\textbf{z}\in\Omega,t>0.
\end{equation}
Since the local human population distribution can reach a time independent dynamic balance in a short time, thus the above parabolic equation degenerates to  the following elliptic equation
\begin{equation}\label{human_distr01}
\Delta x_3+x_3=0,\textbf{z}\in\Omega.
\end{equation}
Moreover, considering the geographical location affect of sea oriented direction we assume that human population $x_3(t,z_1,z_2)=x_3(t,z_1)$ is independent with the coast line direction. And the parameters $d_i, r_i, k_i, h_i(i=1,2.)$ are positive real numbers.

The systems (\ref{prey_predator systems}) can be rewritten as the following PDAEs
\begin{equation}\label{sect5_01}
  E\frac{\partial {x}}{\partial t}=D\Delta{x}+{f}({x})
\end{equation}
where
$$E=\left(
           \begin{array}{ccc}
             1 & 0 & 0 \\
             0 & 1 & 0 \\
             0 & 0 & 0 \\
           \end{array}
         \right),
{x}=(x_1,x_2,x_3)^T,
D=\textrm{diag}(d_1,d_2,1),$$
$${f}({x})=\left(
                    \begin{array}{c}
                      r_1x_1(1-\frac{x_1}{N_1}-k_1 \frac{x_2}{N_2}-h_1x_3) \\
                      r_2x_2(-1+k_2\frac{x_1}{N_1}-\frac{x_2}{N_2}-h_2x_3) \\
                      x_3 \\
                    \end{array}
                  \right).
$$
It is obvious that the positive equilibrium of system (\ref{prey_predator systems}) is the positive solution of the following nonlinear equations
$$1-\frac{x_1}{N_1}-k_1 \frac{x_2}{N_2}-h_1x_3=0, -1+k_2\frac{x_1}{N_1}-\frac{x_2}{N_2}-h_2x_3=0,x_3=0.$$
The positive equilibrium is
$$\left(\frac{N_1(k_1+1)}{1+k_1k_2},\frac{N_2(k_2-1)}{1+k_1k_2},0\right),$$
where $k_2>1$.

\begin{figure}
   \includegraphics[width=0.9\textwidth]{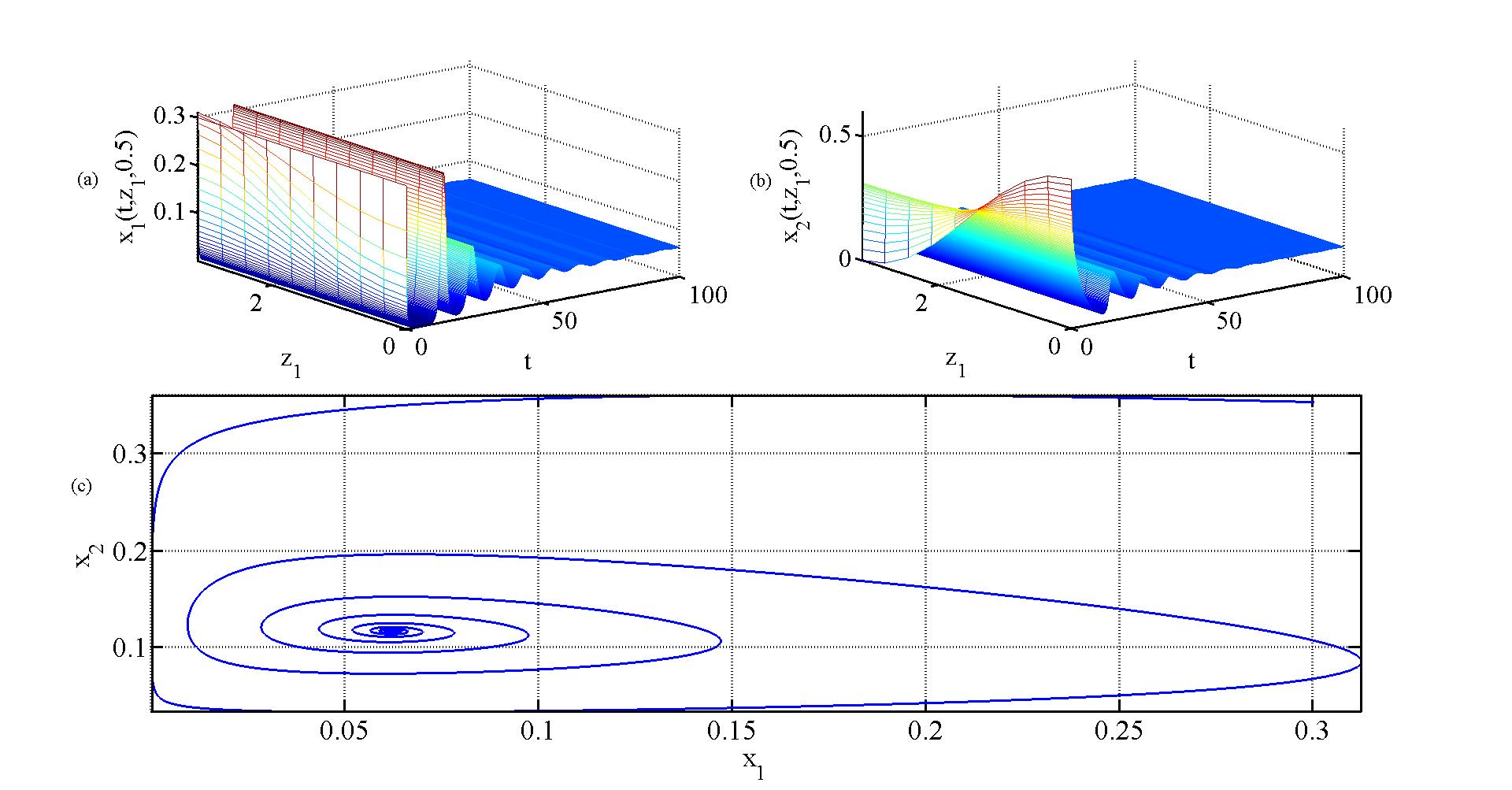}\\
   \caption{\label{simul_result01}Numerical results for the state estimation analysis. (a) Time spatial response $x_1(t,z_1,0.5)$ in the $(t,z_1)$ domain; (b) Time spatial response $x_2(t,z_1,0.5)$ in the $(t,z_1)$ domain; (c)The phase trajectory about $x_1,x_2(h_1=0.2,h_2=0.1).$}
 \end{figure}

\begin{figure}
  \includegraphics[width=0.9\textwidth]{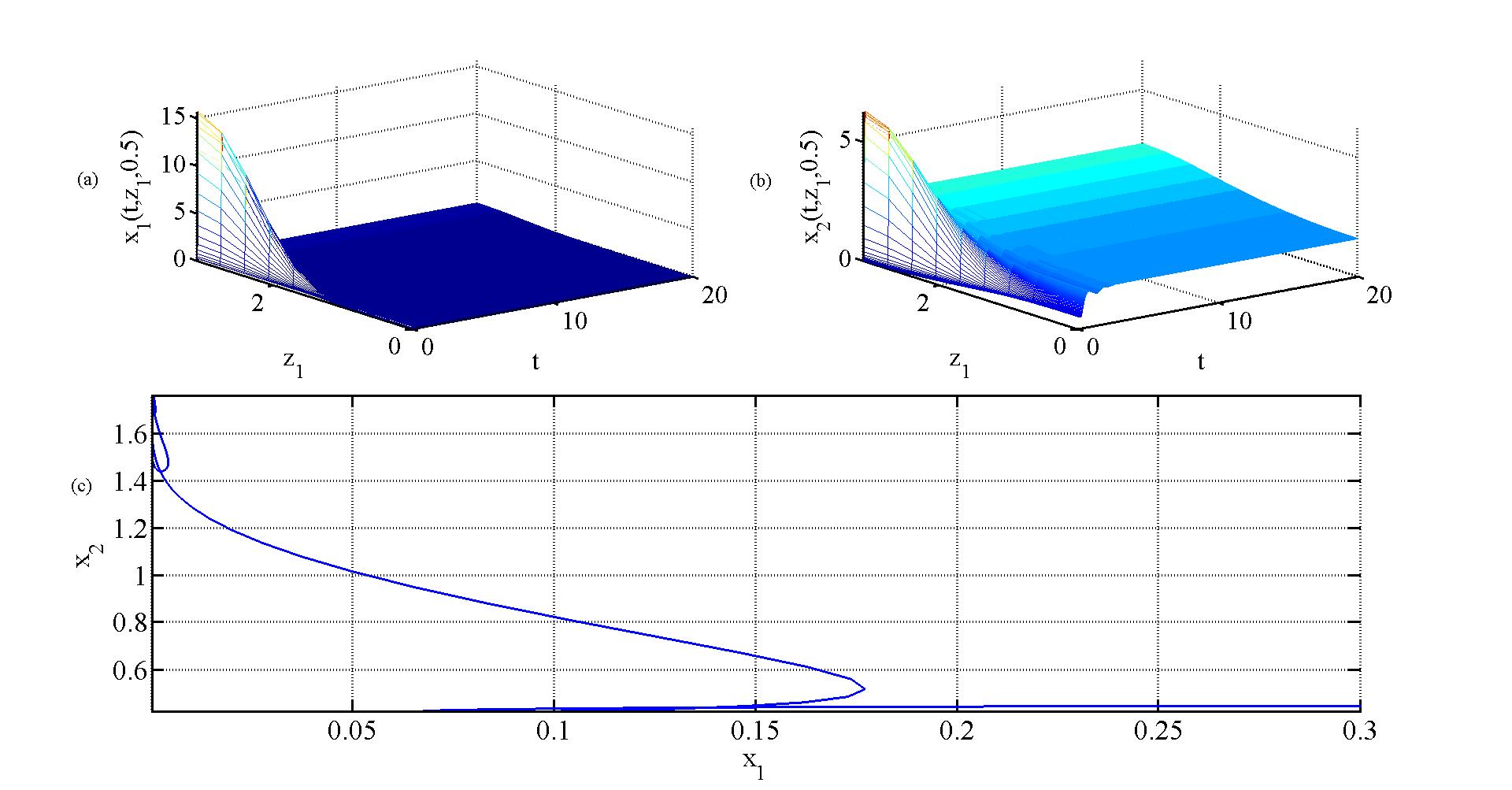}\\
  \caption{\label{simul_result02}Numerical results for the state estimation analysis.(a) Time spatial response $x_1(t,z_1,0.5)$ in the $(t,z_1)$ domain; (b) Time spatial response $x_2(t,z_1,0.5)$ in the $(t,z_1)$ domain; (c)The phase trajectory about $x_1,x_2(h_1=24,h_2=0.1$).}
 \end{figure}

\begin{figure}[htpb]
 \centering
   \includegraphics[width=0.45\textwidth]{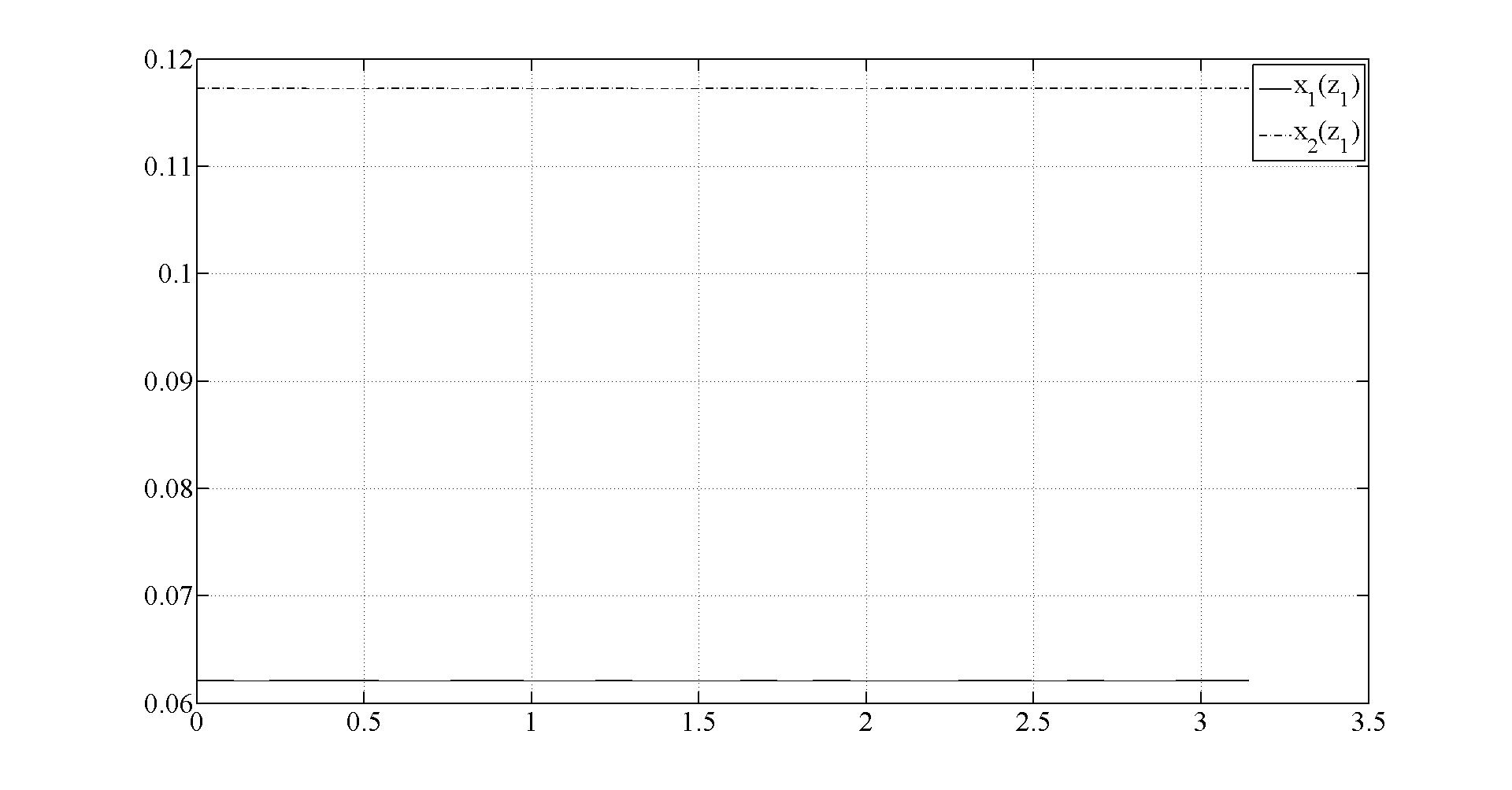}
   \includegraphics[width=0.45\textwidth]{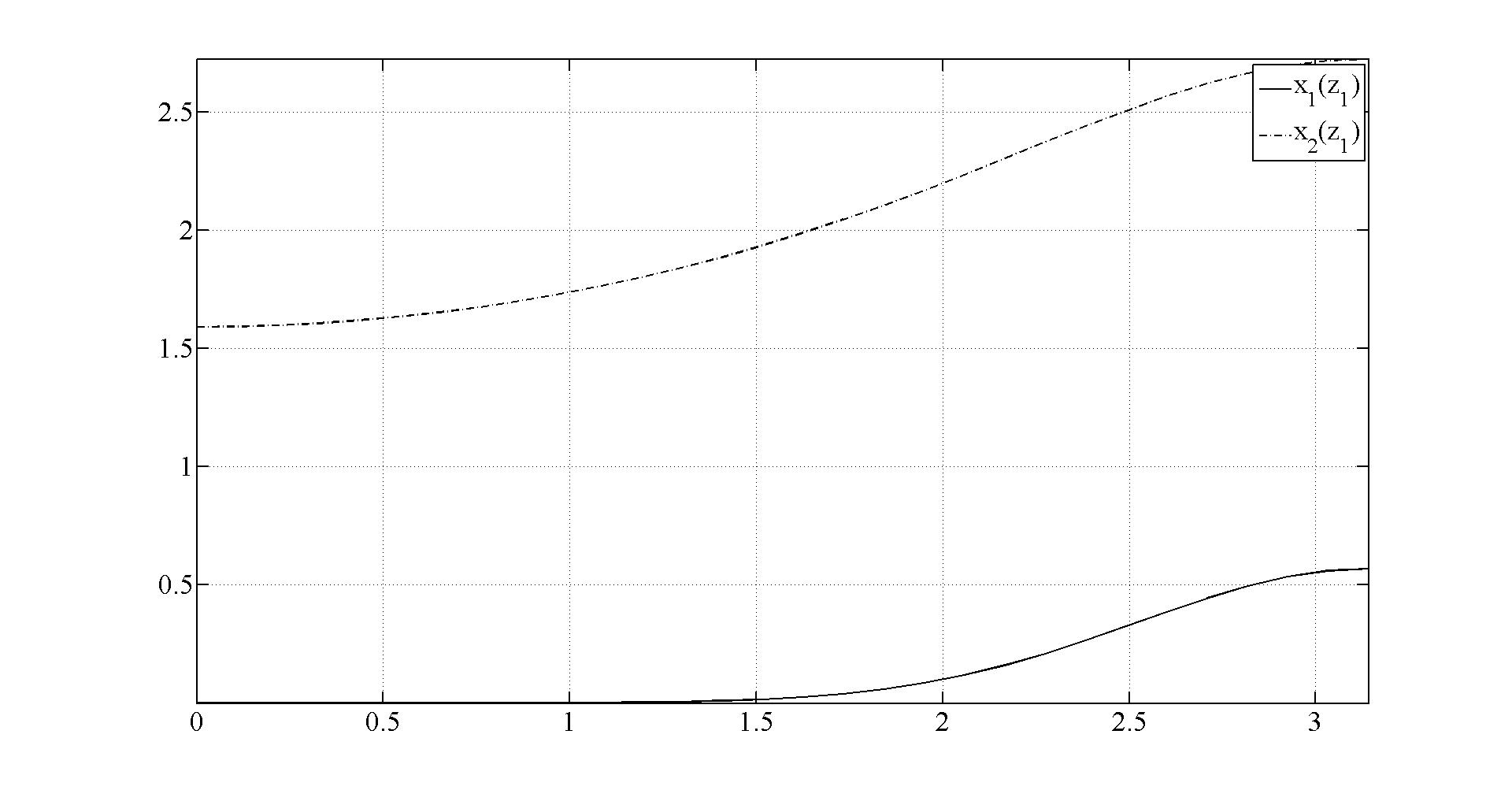}\\
   \caption{\label{simul_result03}Numerical results on spatial convergence of $x_1,x_2.$}
\end{figure}

Now we consider the local stable property of the equilibrium.
We choose $\Omega:=\{(z_1,z_2)|0<z_1<\pi,0<z_2<1\}$. To the consistency conditions the initial conditions are taken as $x_1^0(\textbf{z})=0.3,x_2^0(\textbf{z})=0.3(1+cos(z_1))$. The remaining system parameters are chosen as $N_1=N_2=1$,$r_1=2,r_2=0.2,k_1=8 ,k_2=18,,D=diag(2,3,1)$.
By directly computing we have the Jacobi matrix of the linearized system at the equilibrium
$$A_J=\left(
      \begin{array}{ccc}
        -18/145 & -144/145 & -18h_1/145 \\
        306/725 & -17/725 & -17h_2/725 \\
        0 & 0 & 1 \\
      \end{array}
    \right).
$$
According to the stability analysis of theorem \ref{thm_stable_analy} in section 4, we consider the local stability property of this PDAEs at the equilibrium in different parameter values(see Tab.\ref{data_tabel}).
Due to the particular choice of the system parameters $h_1, h_2$ the state variables $x_1,x_2$ show some different dynamical properties. By increasing the value of $h_1$, the value of $nd_1\mu_1-\|A\|$ can be changed from positive to negative.
Simulation results for the system (\ref{prey_predator systems}) in the $(t,z_1)$ domain at $z_2=0.5$ are depicted in Fig.\ref{simul_result01} and Fig.\ref{simul_result02}. In Fig.\ref{simul_result01} (a) and (b) the state variables $x_1,x_2$ spatial-temporal show the corresponding exponentially stable. It shows great difference with the case in Fig.(\ref{simul_result02}). Fig.(\ref{simul_result03}) shows the different spatial convergence property of the system. Ecologically the stable prey-predator relationship stands when the diffusion rates of prey and predator $diag(d_1,d_2,1)$ are very high with low human affect rates $h_1,h_2$ which correspond with the experience.
\begin{center}
    \begin{table}
    \caption{\label{data_tabel} Values of the simulation parameters}
        \begin{tabular}{c|c|c}
            \hline
            \textbf{Case}& $\|A_J\|$ & $d_1\mu_1-\|A_J\|$ \\   \hline
            $h_1=h_2=0.1$ & 1.0069 & 1.9931 \\  \hline
                $h_1=24,h_2=0.1$ &  3.2841 & -0.7159 \\
            \hline
        \end{tabular}
    \end{table}
\end{center}
\section*{Conclusion}
In this study we have studied the problem of some paraboli-elliptic type PDAEs in high dimensional domain. With the decomposition ideas derived from PDE theory we built the IDAEs to reconstruction of PDAEs. Some spectrum theory result of the IDAEs corresponding the PDAEs are proposed. A exponential stable result on the PDAEs is presented  through the energy estimation about the state variables under the homogenous Neumann boundary conditions for the positive diffusion matrix. Finally as an application we built some wetland conservation model with social behaviour. The numerical results show the effectiveness of this development.
\section*{Conflict of Interests}
The authors declare that there is no conflict of interests
regarding the publication of this paper.
\section*{Acknowledgments}
The research is supported by N.N.S.F. of China
under Grant No. 61273008 and No. 61104003. The
research is also supported by the Key Laboratory of
Integrated Automation of Process Industry
(Northeastern University).The author is grateful to
the anonymous referee for a careful checking of the
details and for helpful comments that allow us to
improve the manuscript.
\bibliographystyle{unsrt}
\bibliography{mybib}
\end{document}